\documentclass[a4paper,12pt]{amsart}
\usepackage[foot]{amsaddr}
\usepackage{graphicx}        % standard LaTeX graphics tool
\usepackage{multicol}        % used for the two-column index
\usepackage[bottom]{footmisc}% places footnotes at page bottom
\usepackage{bm}
\usepackage{graphicx}        % standard LaTeX graphics tool
\usepackage{array}
\newcolumntype{P}[1]{>{\centering\arraybackslash}p{#1}}
\usepackage{array,multirow}
\newcolumntype{C}[1]{>{\centering\arraybackslash}m{#1}}
\usepackage{float}

\usepackage[margin=2.7cm]{geometry}
\usepackage{systeme}
\usepackage{graphicx}
\usepackage[]{units}
\usepackage{color}
\usepackage{mathrsfs}
\usepackage[ansinew]{inputenc}
 \usepackage[breaklinks]{hyperref}
\usepackage{float}
\usepackage{amsmath}
\usepackage{isomath}
\usepackage{multirow}
\usepackage[dvipsnames]{xcolor}
\usepackage{color}

\usepackage{hyperref}

\newcommand{\green}{\color{black}}

\makeindex             % used for the subject index
\begin{document}

\title[]{A Reduced Order Approach for the Embedded Shifted Boundary FEM and a Heat Exchange System on Parametrized Geometries
}

\author{E.N. Karatzas\textsuperscript{1,*}}
\address{\textsuperscript{1}SISSA, International School for Advanced Studies, Mathematics Area, mathLab Trieste, Italy.}

\address{\textsuperscript{2}Civil and Environmental Engineering, Duke University, Durham, NC 27708, United States.}

\thanks{\textsuperscript{*}Corresponding Author.}

\email{efthymios.karatzas@sissa.it}

\author{G. Stabile\textsuperscript{1}}
\email{gstabile@sissa.it}

\author{N. Atallah\textsuperscript{2}}
\email{nabil.atallah@duke.edu}

\author{G. Scovazzi\textsuperscript{2}}
\email{guglielmo.scovazzi@duke.edu}

\author{G. Rozza\textsuperscript{1}}
\email{grozza@sissa.it}

%  % \titlerunning{POD Reduced Basis and SBM for a Heat Exchange System} 
% \author{E.N. Karatzas, G. Stabile, N. Atallah, G. Scovazzi and G. Rozza}

% \institute{Efthimios N. Karatzas \at SISSA, International School for Advanced Studies, Mathematics Area, mathLab, Trieste, 34136, Italy,
%  \email{efthymios.karatzas@sissa.it}
% %
% \and Giovanni Stabile \at SISSA, International School for Advanced Studies, Mathematics Area, mathLab, Trieste, 34136, Italy, \email{gstabile@sissa.it}
% %\and ...?... %L. Nouveau  \at Civil and Environmental Engineering, Duke University, Durham, NC 27708, United States, \email{leo.nouveau@duke.edu}
% %
% \and Nabil Atallah 
% \at Civil and Environmental Engineering, Duke University, Durham, NC 27708, United States, \email{nabil.atallah@duke.edu}
% %
% \and Guglielmo Scovazzi 
% \at Civil and Environmental Engineering, Duke University, Durham, NC 27708, United States, \email{guglielmo.scovazzi@duke.edu}
% \and Gianluigi Rozza \at SISSA, International School for Advanced Studies, Mathematics Area, mathLab, Trieste, 34136, Italy, \email{grozza@sissa.it}
% }

% \maketitle


\date{}

\dedicatory{}

\maketitle

\begin{abstract}
A model order reduction technique is combined with an embedded boundary finite element method with a POD-Galerkin strategy. The proposed methodology is applied to parametrized heat transfer problems and we rely on a sufficiently refined shape-regular background mesh to account for parametrized geometries. In particular, the employed embedded boundary element method is the Shifted Boundary Method (SBM), recently  proposed in \cite{MaSco17_2}. This approach is based on the idea of shifting the location of true boundary conditions to a surrogate boundary, with the goal of avoiding cut cells near the boundary of the computational domain.  
This combination of methodologies has multiple advantages. In the first place, since the Shifted Boundary Method always relies on the same background mesh, there is no need to update the discretized parametric domain. Secondly, we avoid the treatment of cut cell elements, which usually need particular attention. Thirdly, since the whole background mesh is considered in the reduced basis construction, the SBM allows for a smooth transition of the reduced modes across the immersed domain boundary.
The performances of the method are verified in two dimensional heat transfer numerical examples. 
\end{abstract}

\section{Introduction}
In this work we present a reduced order modeling strategy for parametrized geometries, starting from an embedded boundary method solver. The main idea in the current manuscript is to exploit the advantages of embedded methods and in particular of the Shifted Boundary Method (SBM), \cite{MaSco17_2,MaSco17_3,SoMaScoRi17}, in a reduced order modeling setting. Embedded methods, as full order conformal finite element methods, discretize the original set of equations into a usually high dimensional system of algebraic equations. When a large number of different system configurations need to be tested, or a large reduction in computational cost is the goal, the resolution of such high dimensional system of equations becomes unfeasible. Reduced Order Methods (ROM) have demonstrated to be a viable way to limit the computational burden \cite{HeRoSta16,quarteroniRB2016,ChinestaEnc2017,BeOhPaRoUr17}. In this particular case, the attention is focused on parametrized geometries. The methodology is tested on a simple heat transfer problem which will serve as a base for future more complex scenarios such as flow problems \cite{KaStaNoScoRo18}. The manuscript is organized as follows: in Section~\ref{sec:HF} we introduce the mathematical problem and its full order discretization; in Section~\ref{sec:ROM} we present the reduced order model formulation and its main features and differences with respect to a standard setting; finally in Section~\ref{sec:num_exp} numerical results are reported, and in Section~\ref{sec:concl} conclusions and perspectives for future improvements are given.

\section{Full Order Model approximation} \label{sec:HF}
We start recalling, by a sketch description, the continuous strong formulation of the problem and the weak formulation used for the full-order discretizaton of the problems under consideration. The discrete SBM formulation will be used for the Full Order Method (FOM) simulation during the offline stage. The ROM is constructed using a Proper Orthogonal Decomposition (POD) Galerkin approach following what is reported in Section~\ref{sec:ROM}. 

\subsection{The Thermal-Heat exchange model}
Given a $k-$dimensional parameter space $\mathcal P$ and the parameter vector $\mu\in \mathcal P \subset \mathbb{R} ^k
$,  let ${\mathcal{D}}(\mu)\subset {\mathbb R}^d$, $d=2,3$ be a bounded parametrized domain depending on $\mu$, with boundary $\Gamma(\mu)$. We consider the following model problem in $\mathcal D(\mu)$: 

{Find the temperature 
%$T(\mu)\in H^2({\mathcal D}(\mu))$ %
$T(\mu):{\bar{\mathcal D}(\mu)}\times \mathcal P \to \mathbb{R}^d$ 
such that in $\mathcal P$ we have }
\begin{eqnarray} \label{P:Poisson}
-\Delta  T(\mu)&=&f(\mu)\quad\,\,\,\, 
\mbox{ in }
\mathcal{D}(\mu),
\nonumber\\
T(\mu)&=&g_D(\mu) \quad   \mbox{ on }  \Gamma_D(\mu),
\end{eqnarray}
where $\Gamma_D(\mu)$ is the boundary onto which a Dirichlet boundary condition is applied, and the imposed forces $f(\mu)$, $g_D(\mu)$ are given functions in $\mathcal D(\mu)$ and on the boundary $\Gamma_D(\mu)$, respectively.

\subsection{%Continuous 
Weak SBM formulation}
\label{sec:2}
In this subsection we briefly recall the SBM formulation which was originally presented in \cite{MaSco17_2, MaSco17_3, SoMaScoRi17}. In what follows, we denote by $\tilde  \Gamma$ the surrogate boundary composed of the edges/faces of the mesh that are the closest to the true boundary $\Gamma$. The closest faces/edges of $\tilde \Gamma$ to $\Gamma$ are detected  by means of a closest-point projection algorithm. 

\begin{figure}[t]
\centering
\begin{minipage}{0.9\textwidth}
\centering
\includegraphics[width=0.25\textwidth]{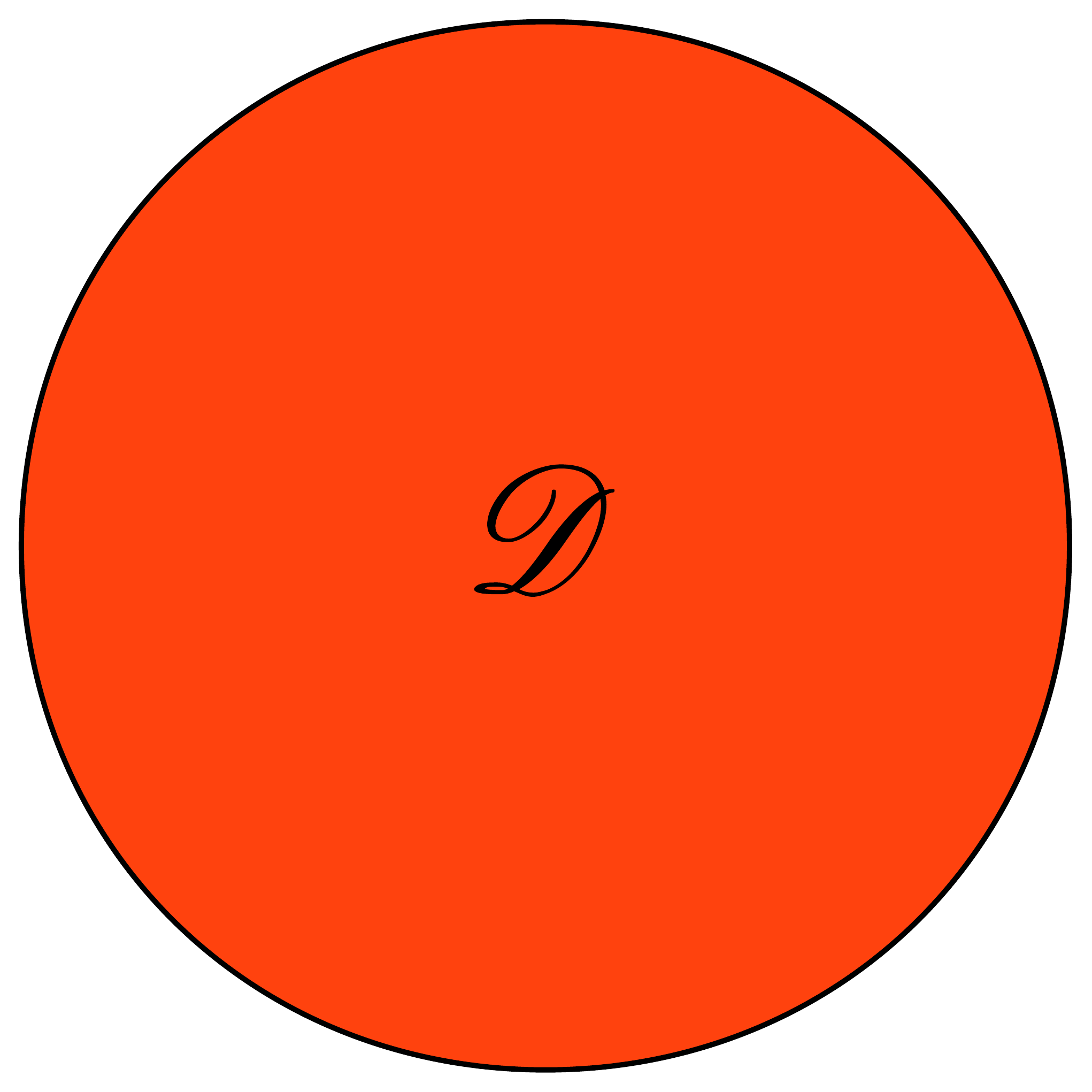}
\includegraphics[width=0.25\textwidth]{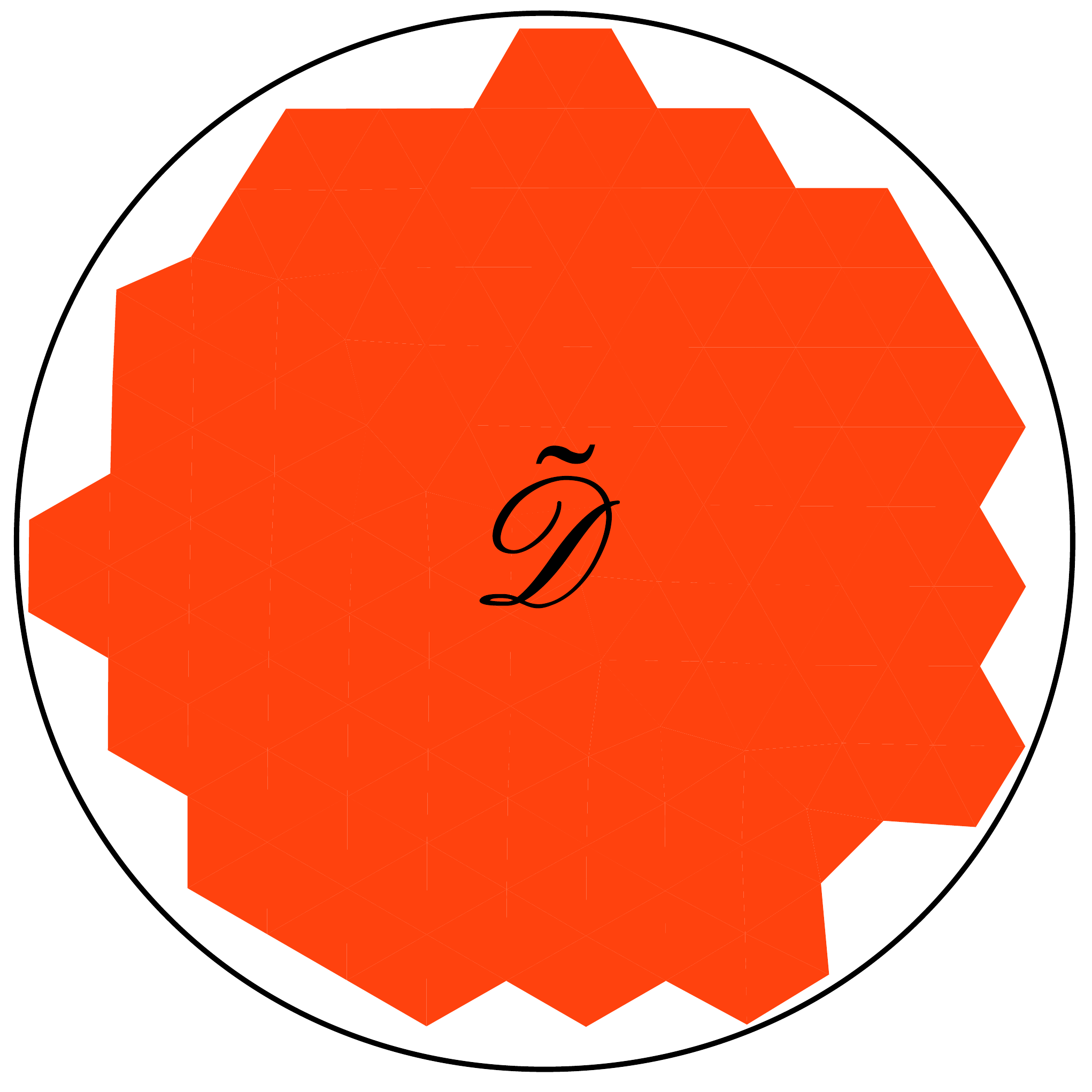}
\includegraphics[width=0.25\textwidth]{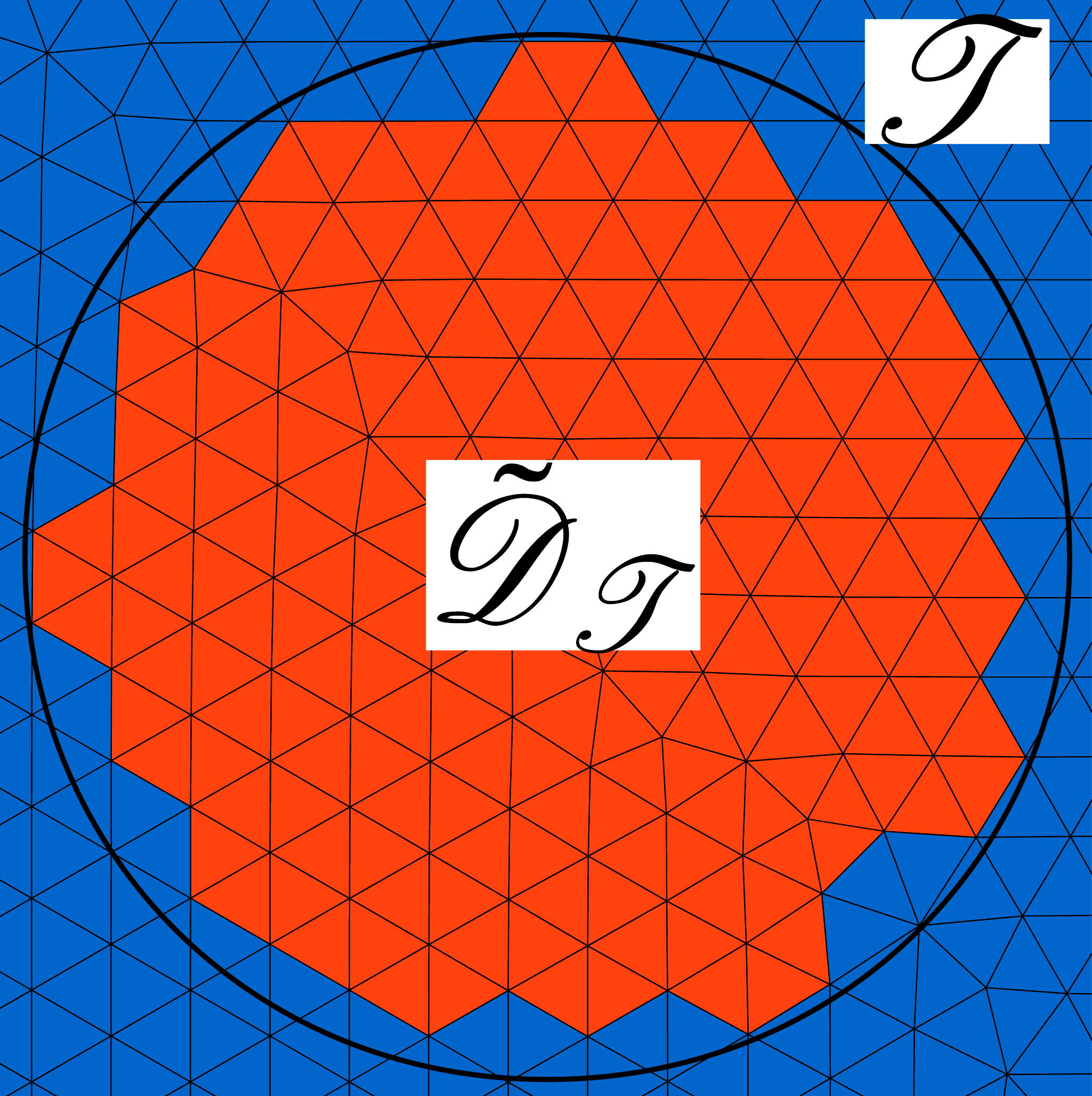}\\
\includegraphics[width=0.4\textwidth]{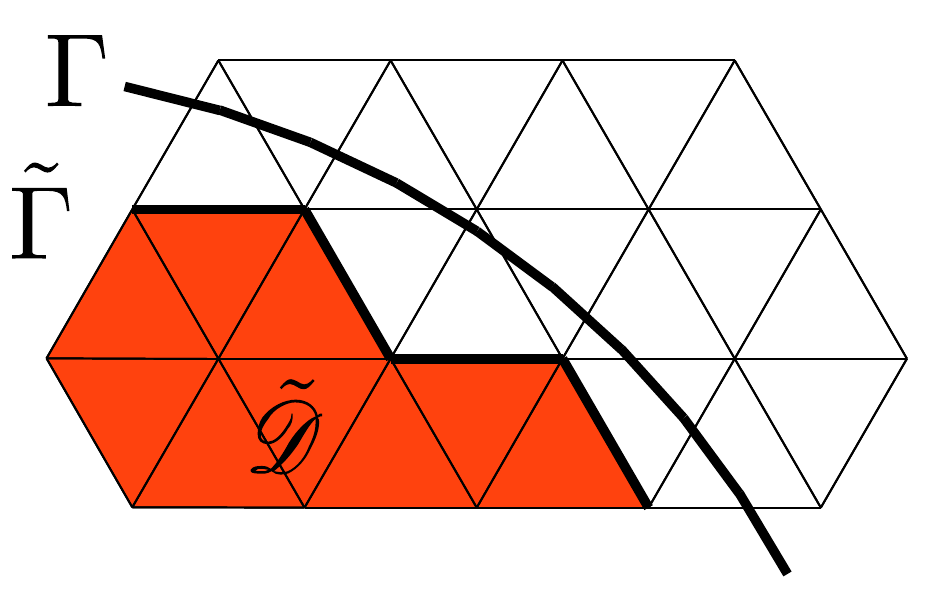}
\includegraphics[width=0.4\textwidth]{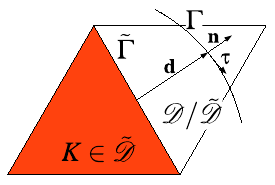}
\end{minipage}
\caption{Example of the SBM mesh on a disc. In the first row, from left to the right: the real geometry; the SBM surrogate geometry and the background mesh together with the surrogate SBM discretized geometry. In the second row, from left to the right: a zoom of the surrogate SBM mesh/ surrogate boundary and the normal and distance vector considering one element.}
\label{fig:SurrogateMesh}       % Give a unique label
\end{figure}

The surrogate boundary $\tilde \Gamma$ encloses the surrogate domain $\tilde {\mathcal{D}}$. Furthermore, $\bf{\tilde{n}}$ indicates the unit outward-pointing normal to the surrogate boundary $\tilde \Gamma$, and it differs from the outward-pointing normal  $\bf  n$  of ${\Gamma}$ (see Figure \ref{fig:SurrogateMesh}). 

Notice also that the closest-point projection, in spite of the segmented/faceted nature of the surrogate boundary $\tilde \Gamma$ is actually a smooth mapping $\bf M$ from points $\tilde{x}$ on $\tilde \Gamma$ to points $x$ on $ \Gamma$, namely,
\begin{eqnarray*}
{\bf{M}}: {\bf \tilde{x}}|_{\tilde\Gamma}\to {\bf x}|_{\Gamma},
\end{eqnarray*}
which also defines a distance vector function:
$$\bf d \equiv d_M (\tilde{x}) = x - \tilde x = [M - I]( \tilde x). $$
The distance vector, as seen in Figure \ref{fig:SurrogateMesh}, is oriented along the normal to the true boundary, that is $\bf d = ||\bf d||\bf n$, as a consequence of the use of the closest point projection. Between the normal $\bf n$ to the true boundary and the normal $ {\bf {\tilde n}}$ to the surrogate boundary, the minimal grid resolution assumption $\bf n\cdot \tilde n \ge 0$ is made. The unit normal vector $\bf n$ and the unit tangential vectors $\bm\tau_i$   $(1 < i < d - 1)$ to the boundary $\Gamma$, can be easily extended to the boundary $\tilde \Gamma$ since
$\bf \bar n(\tilde x) \equiv n(M(\tilde x))$,
 {{\green{${\bar{\bm\tau} _i \bf({{\tilde x}})} \equiv{\bm \tau}_i (\bf M(\tilde x))$}}. Here we denote by $\bf \bar n$, $\bf\bar{\bm{\tau}} _i$ the extensions to $\tilde{\Gamma}$ of $\bf n$, ${{\bm\tau} _i}$, which are defined on $\Gamma$.
 In the following, whenever we write $\bf n(\tilde x)$ we actually mean $\bf \bar n(\tilde x)$ at a point ${\bf \tilde x}\in \tilde\Gamma$,
 and similarly for $ {\bm\tau} _i ( {\bf\tilde x})$ and $ \bf\bar{\bm\tau}_i ( \bf\tilde x)$. 
Moreover, the above constructions are the key ingredients when building an extension $\bar g_D$ of the Dirichlet boundary condition $g_D$ to the boundary  $\tilde  \Gamma$ of the surrogate domain. 

Now we can introduce the Shifted Boundary (SB) variational formulation. The SBM weak discrete formulation for the heat exchange system, with non-homogeneous Dirichlet boundary conditions, reads: 

Find $T \in 
V_h %(\tilde{\mathcal D}(\mu))
= \left\{ \upsilon \in C^0 (%\bar
{\tilde{\mathcal{D}}}(\mu))^{} :  \upsilon|_K \in {P}^1 (K)^{}, \forall K \in {\tilde{\mathcal{D}}}_{\mathcal T}(\mu) \right\}
 $, %V_h %(\tilde{\mathcal D}(\mu))
%= \left\{ v(\mu) \in C^0 (%\bar
%{\tilde{\mathcal{D}}}(\mu))^{} :  v(\mu)|_K \in {P}^1 (K)^{}, \forall K \in {\tilde{\mathcal{D}}}_{\mathcal T}(\mu) \right\}
% $, 
with number of degrees of freedom equal to $\dim V_h = N_h<\infty$ for all $h>0$ such that %$\forall w\in  V_h% (\tilde{\mathcal D})$
\begin{equation}\label{eq:weak_form}
  a(T, w) = \ell(w),\, \forall w\in  V_h,% (\tilde{\mathcal D})$$
\end{equation} 
with
\begin{eqnarray*}
&&a(T, w) = (\nabla w , \nabla T )_{\tilde{\mathcal D}}
- \langle w + \nabla w \cdot {\bf d}, \nabla T \cdot {\tilde {\bf{n}}}\rangle _{\tilde{\Gamma}_D}
- \langle \nabla w\cdot {\tilde {\bf{n}}}, T + \nabla T\cdot {\bf{d}}\rangle_{ \tilde{\Gamma}_D}
\\
&&\qquad\qquad\,\,\,+ \langle \nabla w \cdot {\bf d}, ({\bf n} \cdot {\tilde {\bf{n}}})/||{\bf{d}}|| \nabla T \cdot {\bf d}\rangle_{\tilde{\Gamma}_D}
+\langle  \alpha / {h^\perp} (w + \nabla w  \cdot {\bf{d}}), T + \nabla T \cdot {\bf d}\rangle _{\tilde{\Gamma}_D},
\\
&&\ell(w) = (w, f)_{\tilde{\mathcal D}}
- \langle \nabla w \cdot {\tilde{\bf n}},{\bar{g}_D}\rangle _{\tilde{\Gamma}_D}
- \langle \nabla w \cdot {\bf d}, (\nabla {\bar g}_D \cdot {\bm \tau}_i) {\bm \tau}_i \cdot {\tilde{\bf n}}\rangle_{\tilde{\Gamma}_D}
\\
&&\qquad\qquad\,\,\,
+ \langle \alpha / {h^\perp}(w + \nabla w \cdot {\bf d}), \bar{g}_D\rangle_{\tilde{\Gamma}_D},
\end{eqnarray*}
where $\alpha$ is the Nitsche penalty parameter, ${h^\perp}$ is a characteristic length of the elements in the direction orthogonal to the boundary and  $\bf d$, $\mathcal T$, ${\tilde{\mathcal{D}}}_{\mathcal T}$ are the distance vector, the background mesh and the discretized surrogate geometry respectively, see e.g Figure \ref{fig:SurrogateMesh}. Finally, the standard notation $(\cdot, \cdot)_{\tilde{\mathcal{D}}}$, $\langle \cdot, \cdot \rangle _{\tilde\Gamma_D}$ have been used for the $L^2(\tilde{\mathcal{D}})$ and $L^2(\tilde\Gamma_D)$ inner products onto the surrogate geometry $\tilde{\mathcal{D}}$ and $\tilde\Gamma_D$, respectively.

The idea of the Shifted Boundary method is to enforce the Dirichlet boundary conditions weakly on the surrogate domain and to modify the value of the boundary conditions to be imposed by means of a second-order accurate Taylor expansion, that is $T + \nabla T\cdot {\bf{d}} \approx \bar{g}_D$, with the purpose of maintaining overall second-order accuracy with a piecewise linear discretization.

The SBM weak formulation %after discretization 
can be transformed in a system of linear equations and rewritten in matrix form:
\begin{equation}\label{eq:system_linear}
\bm{A}(\mu) {\bm{T}}(\mu)   = \bm{F_g}(\mu) \mbox{,}
\end{equation}
where $\bm{A}(\mu)\in {\mathbb R}^{N_h \times N_h}$ corresponds to the bilinear form $a(\cdot,\cdot)$, ${\bm{T}}(\mu)\in {\mathbb{R}}^{N_h \times 1}$ is the vector of the unknowns and
$\bm{F_g}(\mu)\in {\mathbb{R}}^{N_h \times 1}$ corresponds to the linear form $\ell(\cdot)$.
\section{Reduced order method by a POD-Galerkin technique}\label{sec:ROM}
In this section we briefly recall the POD-Galerkin technique used to generate the reduced order model and we highlight its peculiarities with respect to standard approaches. In general a ROM is a simplification of a % full order model 
FOM that preserves essential behavior and dominant effects, for the purpose of reducing solution time or storage capacity. In particular here we employ a projection-based reduced order model which consists of the projection of the governing equations onto the reduced basis space. 

In the recent past, RB methods were applied to linear elliptic equations in \cite{Rozza2008229}, to linear parabolic equations in \cite{grepl2005} and to non-linear problems in \cite{Veroy2003,Grepl2007}. Although the number of works on reduced order models is now considerable (see e.g. \cite{HeRoSta16} and references therein), to the best of the authors' knowledge, only very few research works \cite{BaFa2014} can be found dealing with embedded boundary methods and ROM.

From a reduced order modeling point of view, our aim is to investigate how ROMs are applied to the SBM and, more generally, to embedded boundary methods. Our main interest is to generate ROMs on parametrized geometries. The SBM unfitted/surrogate mesh finite element method is used to apply parametrization and reduced order techniques considering Dirichlet boundary conditions. 

An important objective is also to test the efficiency of a geometrically parametrized reduced order method without the usage of the transformation to reference domains, which can be an important advantage of embedded methods relying on fixed background meshes. 

Before going into the details, we just remind the basics of the reduced basis method. The first step is the generation of a set of full order solutions of the parametrized problem under the  variation of the parameter values. The final goal of RB methods is to approximate any member of this solution set with a low number of basis functions and is based on a two stage procedure, the offline and the online stage, \cite{quarteroniRB2016, RoVe07, Haasdonk2008}. 
\paragraph{\bf Offline stage} In this stage one performs a certain number of full order solves in order to use the solutions for the construction of a low dimensional reduced basis that approximates any member of the solution set to a prescribed accuracy. It is then possible to perform a Galerkin projection of the full order differential operators, describing the governing equations, onto the reduced basis space in order to create a reduced system of equations. This operation usually involves the solution of a possibly large number of high dimensional problems and the manipulation of high-dimensional structures. The required computational cost is high and therefore this operation is usually performed on a high performance system such as a computer cluster.

\paragraph{\bf Online stage} During this stage, that can be performed also on a system with a reduced computational power and storage capacity, the reduced system of equations can be solved for any new value of the input parameters. This offline-online splitting is effective in many scenarios, such as uncertainty quantification, optimization, real-time control, etc,~\cite{BeOhPaRoUr17, ChinestaEnc2017}. 
\subsection{POD}\label{subsec_POD_theory}
In order to generate the reduced basis space, necessary for the projection of the governing equations, one can find in the literature several techniques such as the POD, the Proper Generalized Decomposition (PGD) and the Reduced Basis (RB) with a greedy sampling strategy. For more details about the different strategies, the reader may see \cite{Rozza2008229,Kalashnikova_ROMcomprohtua,
Chinesta2011,Dumon20111387}. We apply here a POD strategy using the method of snapshots \cite{Sirovich1987}. In order to assemble the snapshots matrix, the full-order model is solved for each $\mu \in \mathcal{K}=\{ \mu^1, \dots, \mu^{N_s}\} \subset {\mathbb{R}^k}%\mathcal{P}
$ where $\mathcal{K}$ is a finite dimensional training set of parameters chosen inside the parameter space ${\mathcal P}%\mathcal{P}
$ and $k$ is the size of the vector $\mu$. 
The  number of snapshots is denoted by $N_s$ and the number of degrees of freedom for the discrete full order solution by $N_h$. The snapshots matrix ${\mathcal{S}}$, is then given by $N_s$ full-order snapshots:
\begin{equation}
\bm{\mathcal{S}} = [{T}(\mu^1),\dots,{T}(\mu^{N_s})] \in \mathbb{R}^{N_h\times N_s}.
\end{equation}
Given a general scalar function ${T}:{\mathcal D} \to \mathbb{R}^d$, with a certain number of realizations ${T}_1,\dots, {T}_{N_s}$, and denoting by $({{\cdot},{\cdot}}) _{{\mathcal D }}$ and $||\cdot||_{L^2({\mathcal D })} $ the $L^2({\mathcal{D}})$  inner product and norm  onto the  geometry ${\mathcal{D}}$, the POD problem consists of finding, for each value of the dimension of POD space $N_{POD} = 1,\dots,N_s$, the scalar coefficients $a_1^1,\dots,a_1^{N_s},\dots,a_{N_s}^1,\dots,a_{N_s}^{N_s}$ and functions ${\varphi}_1,\dots,{\varphi}_{N_s}$ that minimize the quantity:
\begin{eqnarray}\label{eq:pod_energy}
E_{N_{POD}} = \sum_{i=1}^{N_s}||{{T}}_i-\sum_{k=1}^{N_{POD}}a_i^k {{\varphi}}_k||^2_{L^2({\mathcal D })}, &&\forall
N_{POD} = 1,\dots,N\\\nonumber
&& \mbox{with } ({{\varphi}_i,{\varphi}_j}) _{{\mathcal D }} = \delta_{ij}, \mbox{\hspace{0.3cm}}  \mbox{}\forall
 i,j = 1,\dots,N_s .
\end{eqnarray}
It can be shown \cite{Kunisch2002492} that the minimization problem of equation~(\ref{eq:pod_energy}) is equivalent of solving the following eigenvalue problem:
\begin{equation}
{\bm{C}}\bm{Q} = \bm{Q}\bm{\lambda} ,\quad \mbox{\hspace{0.5cm} for }{C}_{ij} = ({{T}_i,{T}_j}) _{{\mathcal D }} \mbox{,\, } i,j = 1,\dots,N_s ,\nonumber
\end{equation}
where ${\bm{C}}$ is the correlation matrix obtained starting from the snapshots $\bm{\mathcal{S}}$, $\bm{Q}$ is a square matrix of eigenvectors and $\bm{\lambda}$ is a diagonal matrix of eigenvalues. 

The basis functions can then be obtained with: 
\begin{equation}
{\varphi_i} = \frac{1}{N_s\lambda^{1/2}_{ii}}\sum_{j=1}^{N_s} {T}_j Q_{ij}.
\end{equation}

The POD space are constructed using the aforementioned methodology resulting in the space:
\begin{eqnarray}
\bm{L} = [{{\varphi}}_1, \dots , {{\varphi}}_{N^r}] \in \mathbb{R}^{N_h \times N^r},
\end{eqnarray}
where $N^r < N_s$ is chosen according to the eigenvalue decay of $\bm{\lambda}$, see for example \cite{Rozza2008229,BeOhPaRoUr17}.
\subsection{Main differences with respect to a reference domain approach}
We highlight here that using an embedded approach there is no need to map all the parametrized geometries to a common reference domain as usually done in the reduced order modeling community \cite{%HeRoSta16,
RoVe07,Rozza2009,ballarin2015supremizer,RoHuMa13,Rozza2008229,BeOhPaRoUr17}. The linear and bilinear forms of equation~(\ref{eq:weak_form}), rewritten in a reference domain setting and in a conformal classical finite element method formulation with homogeneous Dirichlet boundary conditions, are transformed into:
\begin{align*}
&\tilde a({w},{T};\mu)=\tilde \ell( w;\mu), \\
&\tilde a({w},{T};\mu) = \int_{{\mathcal{D}^*}} \nabla{w} (J_T(\mu) )^{-1}(J_T(\mu) )^{-T}|J_T(\mu)| \nabla T \mbox{d} {x} , \\
&\tilde \ell(w;\mu) = \int_{\mathcal{D}^*} |J_T(\mu)| {f} {w} \mbox{d} {x},
\end{align*}
where for a reference domain configuration ${\mathcal{D}}^*$, $J_T(\mu)$ and $|J_T(\mu)|$ are the Jacobian of the transformation map $\mathcal{T_M}(\mu): \mathcal{D}^* \to {\mathcal{D}}(\mu)$ 
 and its determinants respectively. For simple geometrical parametrizations, it is possible to find an affine decomposition of the map and therefore of the differential operator ensuring a complete splitting between the offline and the online procedure, see e.g. \cite{HeRoSta16}. For more complex cases such an operation becomes not trivial and therefore, in order to ensure an efficient splitting one has to rely on empirical interpolation techniques or similar methods, \cite{BARRAULT2004667,Carlberg2013623,Rozza2009}. In the proposed method, even though an efficient splitting is not trivial, there is no need to rely on a transformation map. % \cite{}.

All the solutions are in fact referred to a common background mesh and therefore the projection step and the reduced basis generation become straightforward. Each snapshot however has an ``out-of-interest'' region which lives inside the embedded domain and that is usually referred as ``ghost area''. The location of such part of the domain depends on the parameter $\mu$ but the value assumed by the nodes inside that area is arbitrary. The shifted boundary method used herein has the particular advantage that the solution smoothly decreases to zero from the boundary to the interior of the ghost area (see Figure \ref{fig:poisson_zoom}). Besides the closest points, where we have such smooth decrease, the value inside the ghost area is set to zero. Since this choice is arbitrary, other choices are also possible (see \cite{KaBaRo18} for more details). Using such an approach we remark that it is usually not possible to easily recover an affine decomposition of the differential operator with respect to the geometrical parameters. However, as highlighted in the next section, it is still possible to rely on hyper reduction techniques, \cite{Xiao20141,BARRAULT2004667,Carlberg2013623}. 
\begin{figure}[t]
\includegraphics[width=0.40\textwidth]{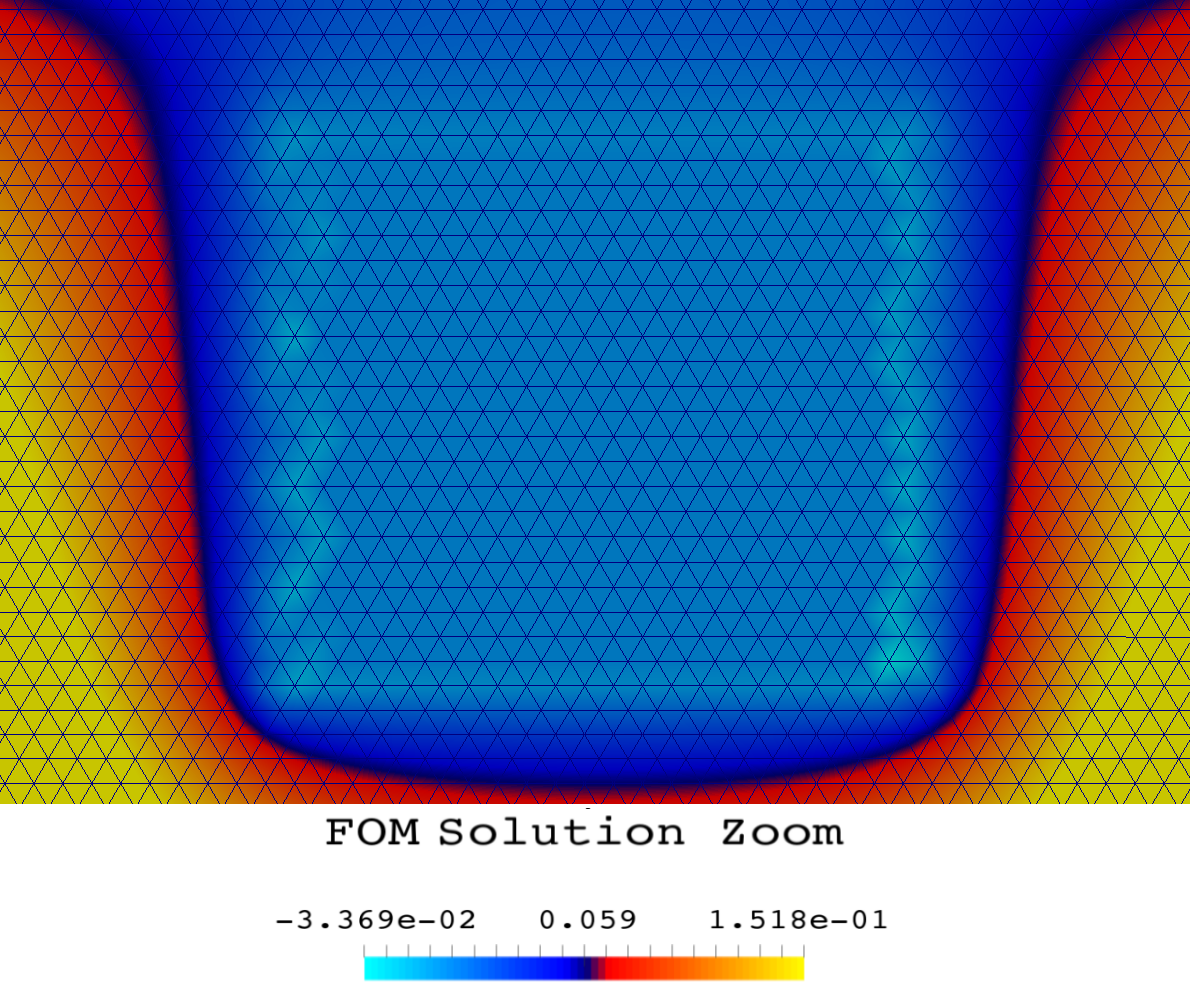}
\caption{A zoom into the embedded rectangle in order to show the smoothing procedure employed by the SBM method inside the ghost area.}
\label{fig:poisson_zoom}
\end{figure}
\subsection{The projection stage and the generation of the ROM}
Once the POD functional space is set, the reduced field can be approximated with: 
\begin{equation}\label{eq:aprox_fields}
{T^r} \approx \sum_{i=1}^{N^r} a_i(\mu) {\varphi}_i(\bm{x}) = {\bm L} \bm{a}(\mu),
\end{equation} 
where the reduced solution vectors $\bm{a} \in \mathbb{R}^{N^r\times1}$ depend only on the parameter values and the basis functions ${\varphi}_i$ depend only on the physical space. The unknown vector of coefficients $\bm{a}$  can then be obtained through a Galerkin projection of the full order system of equations onto the POD reduced basis space and with the resolution of a consequent reduced algebraic system:
\begin{equation}\label{eq:system_linear_proj}
 {\bm L}^T   \bm{A}(\mu)   {\bm L}   {{\bm{a}}(\mu)} = {\bm L}^T \bm{F}(\mu) \mbox{,}
\end{equation}
which leads to the following algebraic reduced system:
\begin{equation}\label{eq:system_linear_reduced}
 \bm{A}^r(\mu) {\green{\bm{a}(\mu)}} = \bm{F}^r (\mu)\mbox{,}
\end{equation}
where $\bm{A}^r(\mu) \in {\mathbb{R}}^{N^r \times N^r}$,  and $\bm{F}^r(\mu) \in {\mathbb{R}}^{N^r \times 1}$ are the reduced discretized operators and reduced forcing vector respectively. The dimension of the reduced operator, as seen also in the numerical examples, is usually much smaller than the dimension of the full order system of equations and therefore much faster to solve. We remark here that the full order discretized differential operators that appear in equation~(\ref{eq:system_linear}) are parameter dependent and therefore, also at the reduced order level, in order to compute the reduced differential operator, we need to assemble the full order operators. Possible ways to avoid such potentially expensive operation, relying on an affine approximation of the full order differential operator, could be to use hyper reduction techniques%, \cite{Xiao20141,BARRAULT2004667,Carlberg2013623}
. In this work, since the attention is mainly devoted to the methodological development of a reduced order method in an embedded boundary setting, rather than in its efficiency, we do not rely on such hyper reduction techniques and we assemble the full order differential operators also during the online stage. Considering that the most demanding computational effort is spent during the resolution of the full order problem rather than in the assembly of the differential operators, as reported in Section~\ref{sec:num_exp}, it is anyway possible to achieve a computational speedup, and the related results are reported in the next section. 

\section{Numerical experiments}\label{sec:num_exp}
{\green{
We consider a parameter space $\mathcal P$ and parameter vector $\mu\in \mathcal P \subset \mathbb{R}$. Let  ${\mathcal{D}}(\mu)\subset {\mathbb R^2}$, be a bounded parametrized domain depending on $\mu$, with boundary $\Gamma_D(\mu)$. In this Section, we report numerical results for the model problem: {Find the reduced basis temperature 
%$T(\mu)\in H^2({\mathcal D}(\mu))$ %
$T(\mu):{\bar{\mathcal D}(\mu)}\times \mathcal P \to \mathbb{R}$ 
such that in $\mathcal P$ we have }
\begin{eqnarray*}
-\Delta  T(\mu)&=&f(\mu)\quad\,\,\,\, 
\mbox{ in }
\mathcal{D}(\mu),
\nonumber\\
T(\mu)&=&g_D(\mu) \quad   \mbox{ on }  \Gamma_D(\mu),
\end{eqnarray*}
where $\Gamma_D(\mu)$ is the embedded boundary onto which a Dirichlet boundary condition is applied, and the imposed forces $f(\mu)$=1, $g_D(\mu)=0$ are forcing data in $\mathcal D(\mu)$ and on $\Gamma_D(\mu)$, respectively. Two different geometries and parameterizations on an embedded rectangle will be examined. In the first example the $y$-coordinate of the embedded domain center is parametrized, and in the second one its aspect ratio is considered as a parameter. 
}}
\subsection{Embedded rectangle with parameterized center}
\begin{figure}[t]
% to insert the figure file.second
\begin{minipage}{\textwidth}
\centering
%\hskip-80pt
%\footnotesize
%(b)
\includegraphics[width=0.51\textwidth]{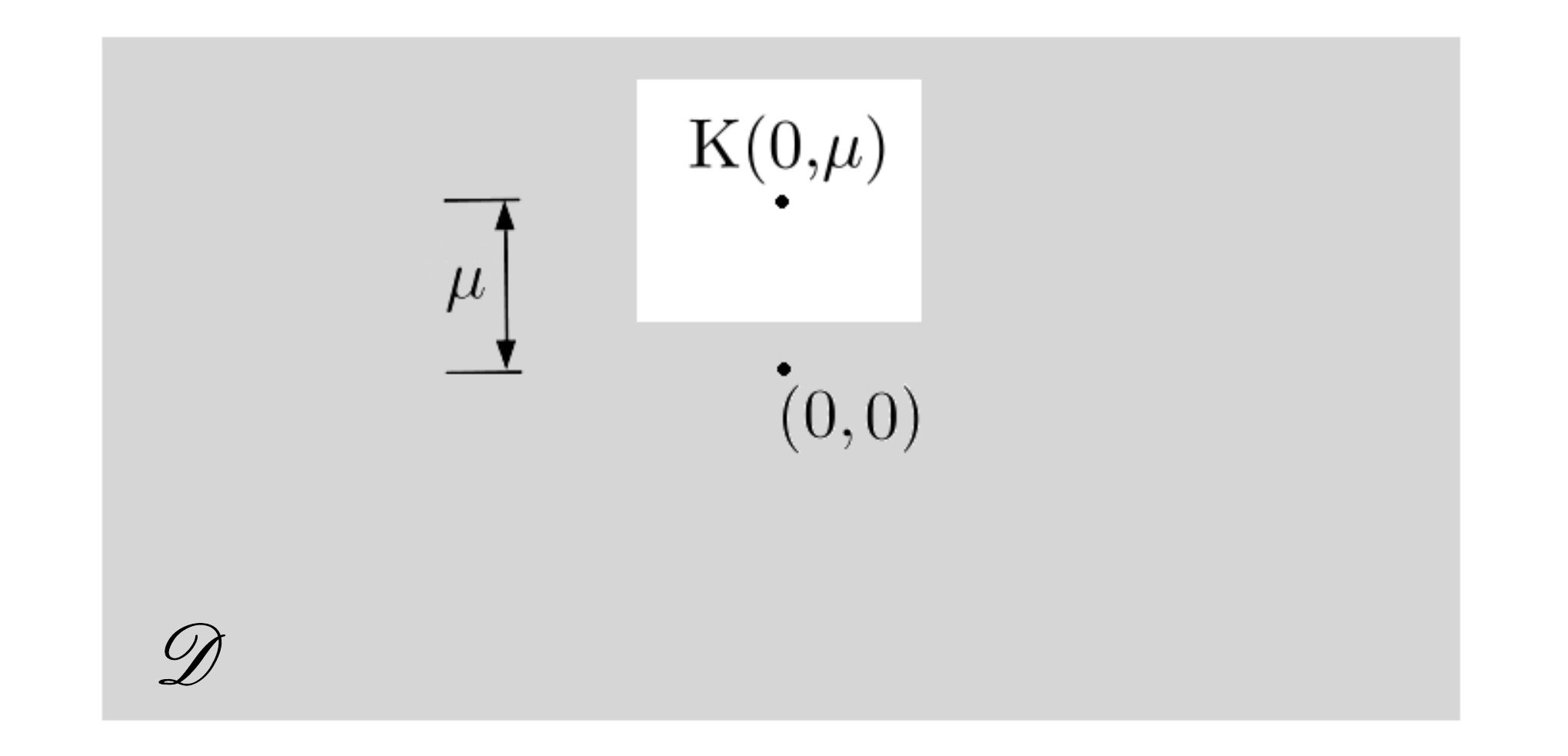}
\end{minipage}
\caption{{\green{Background and embedded geometry for a moving rectangle where the y-coordinate of its baricenter has been parametrized.}}}\label{fig:Truth_Geometry}
\end{figure}
{\green{In this first experiment}} the embedded domain consists of a rectangle of size $0.8\times 0.7$ and its position inside the domain is parametrized with a geometrical parameter $\mu$ which describes the position of the rectangle embedded domain with respect to its y-center {\green{as in Figure \ref{fig:Truth_Geometry}}}. 

The horizontal coordinate of the center of the box is not parametrized and is located in the x-center of the domain. The ROM has been trained with $100$ and $400$ samples for $\mu \in [-0.5,0.5]$ chosen randomly inside the parameter space. To test the accuracy of the ROM we compared its results on $50$ additional samples that were not used to create the ROM and were selected randomly within the same range. The background domain size is {\green{a rectangle of size $[-2,2]\times[-1,1]$ discretized with mesh size $h=0.035$,} while the background mesh boundary is handled as a wall having zero temperture}. 

In Figure \ref{Fig:Modes}, we plot the first four modes obtained with the POD procedure. In Figure \ref{fig:poisson_field}, we plot the full order solution, the reduced solution and the error for the scalar geometrical parametrized heat equation problem and it is possible to notice that the full and the reduced solution are qualitatively indistinguishable. To verify the behavior of the ROM and its sensitivity with respect to the number of modes in Figure \ref{fig:poisson_results} {\green{$\bf{(i)}$}} we compare, for different number of modes, the average of the $L^2$ norm relative error for the $50$ different samples used to test the ROM. The plot is reported for both the simple $L^2$ projection of the full order results on the POD basis functions, and for the ROM results.
\begin{figure}[t]
%\sidecaption[t]
\begin{minipage}{\textwidth}
\includegraphics[width=0.24\textwidth]{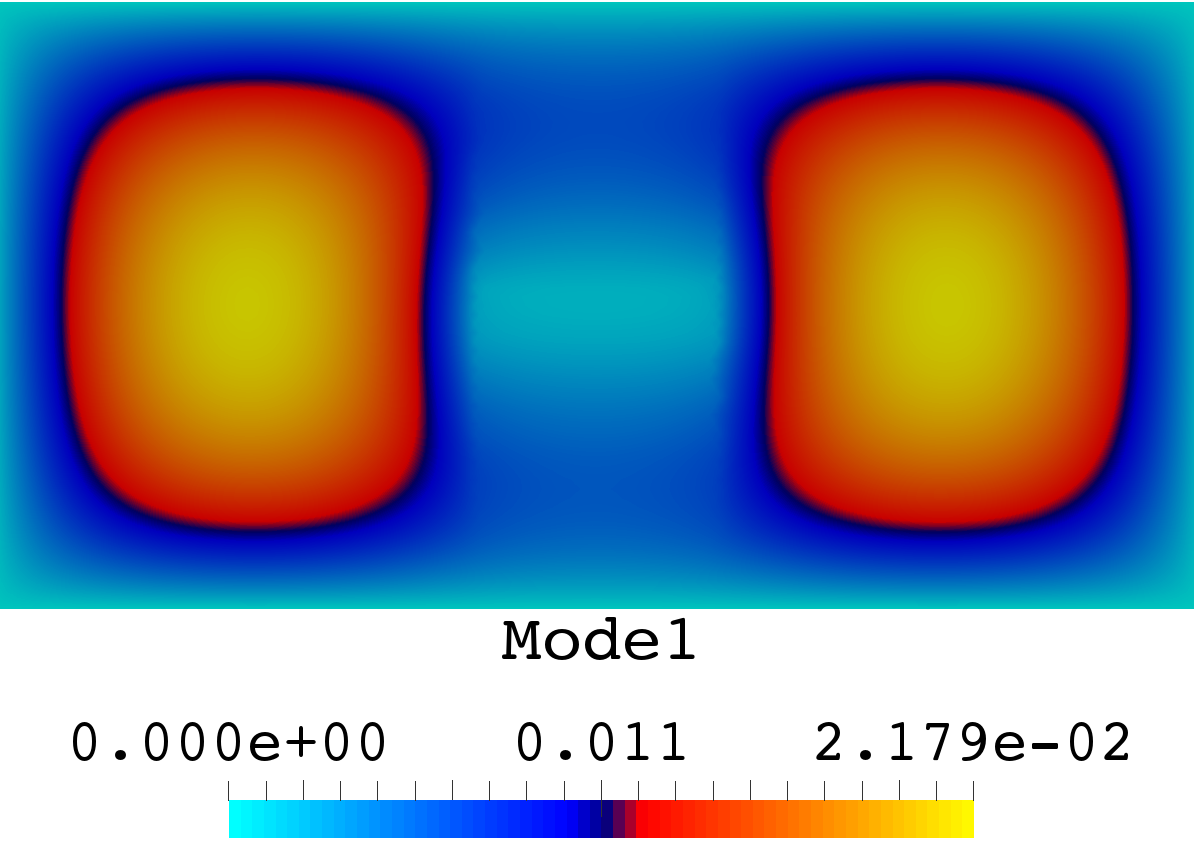}
\includegraphics[width=0.24\textwidth]{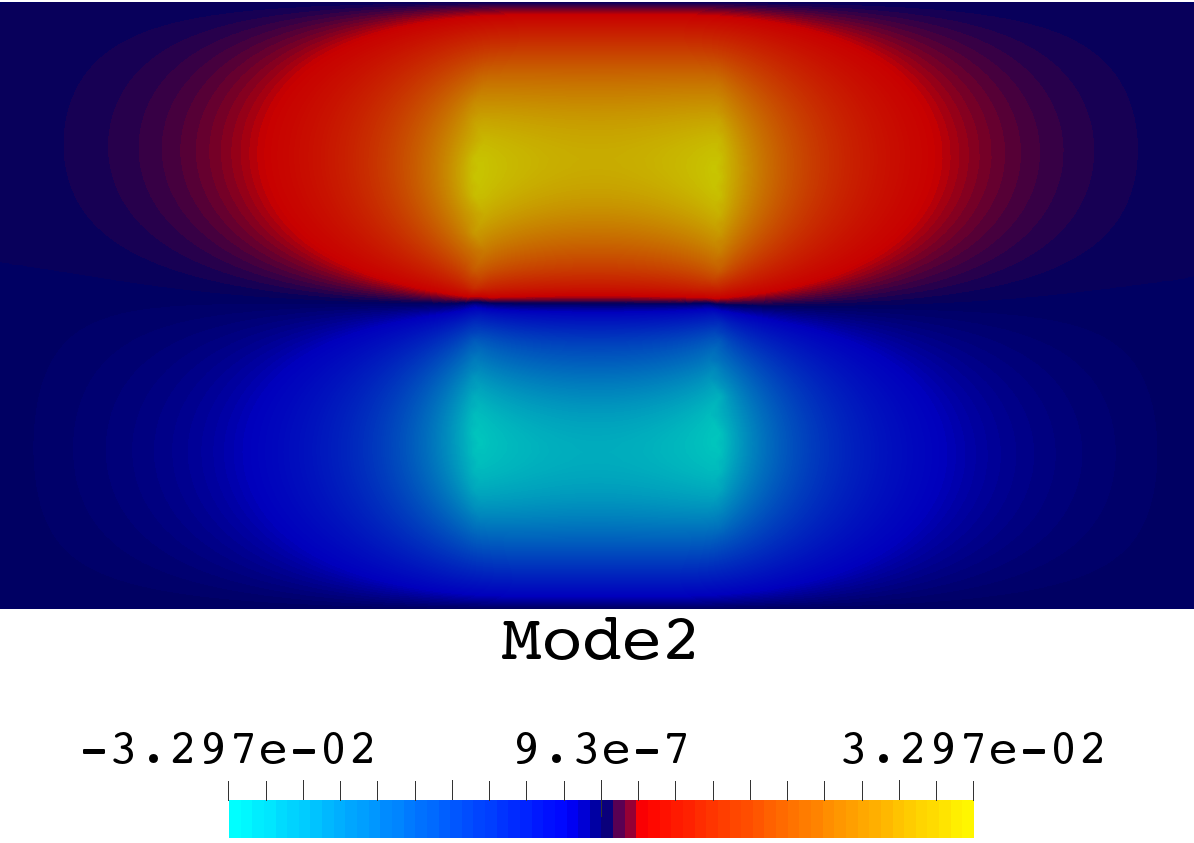}
\includegraphics[width=0.24\textwidth]{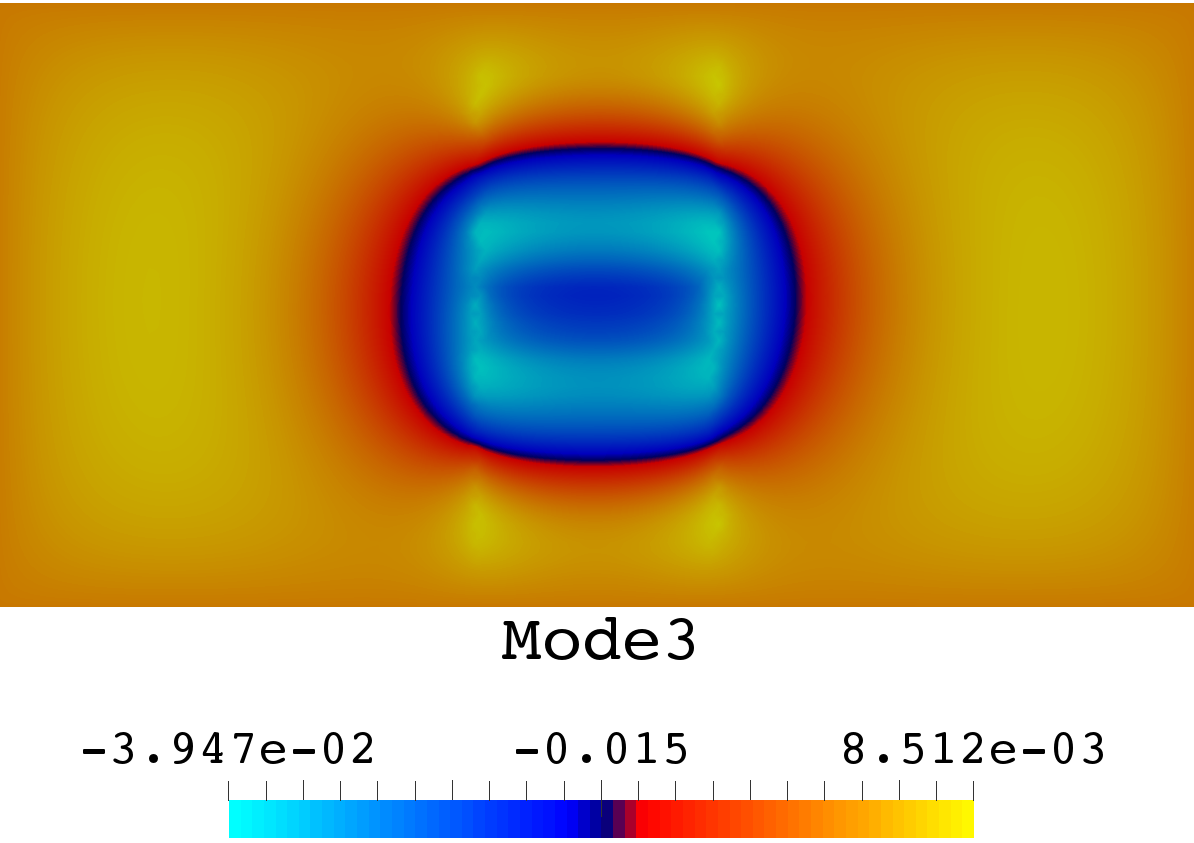}
\includegraphics[width=0.24\textwidth]{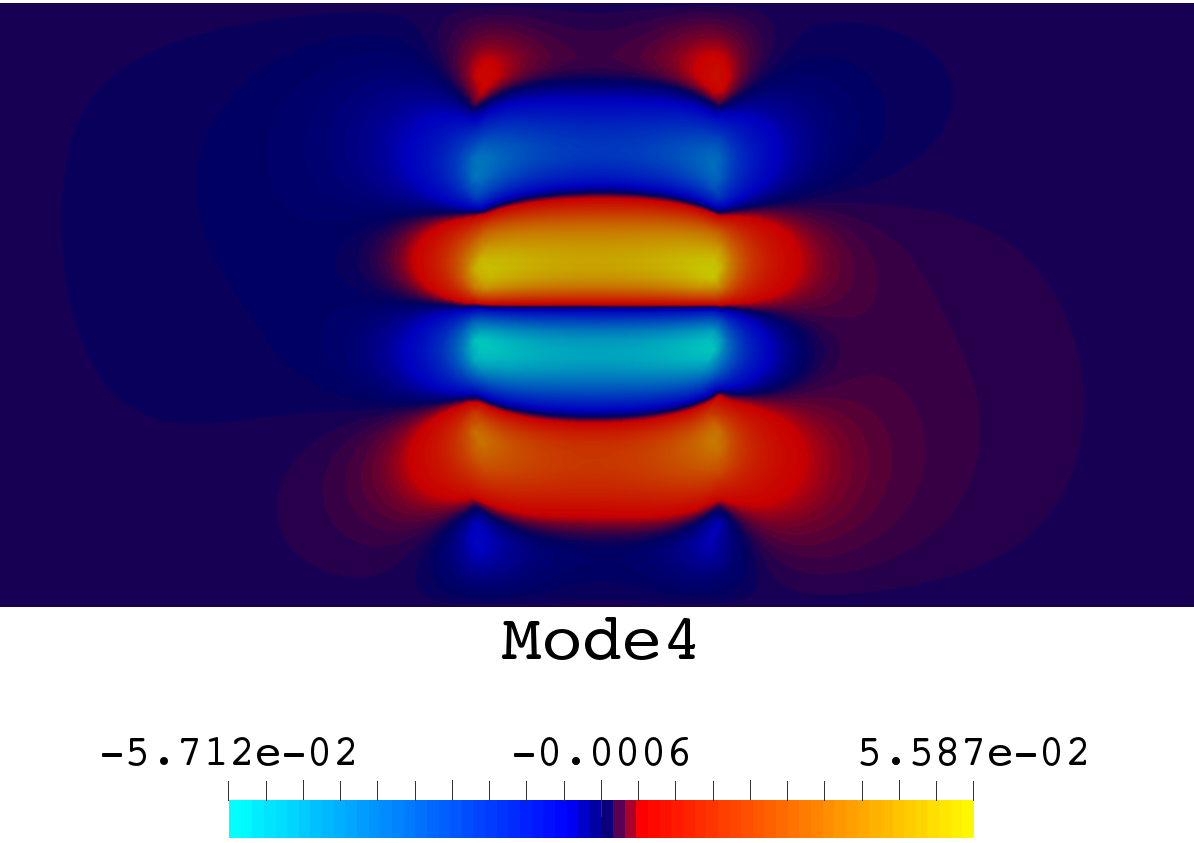}
\end{minipage}
\caption{The first four basis components with $\mu \in [-0.5, 0.5]$ using 100 snapshots in the offline stage.}
 \label{Fig:Modes}
\end{figure}
\begin{figure}[t]
%\sidecaption[t] 
\begin{minipage}{\textwidth}
$\bf{(i)}$ \newline
\includegraphics[width=0.48\textwidth]{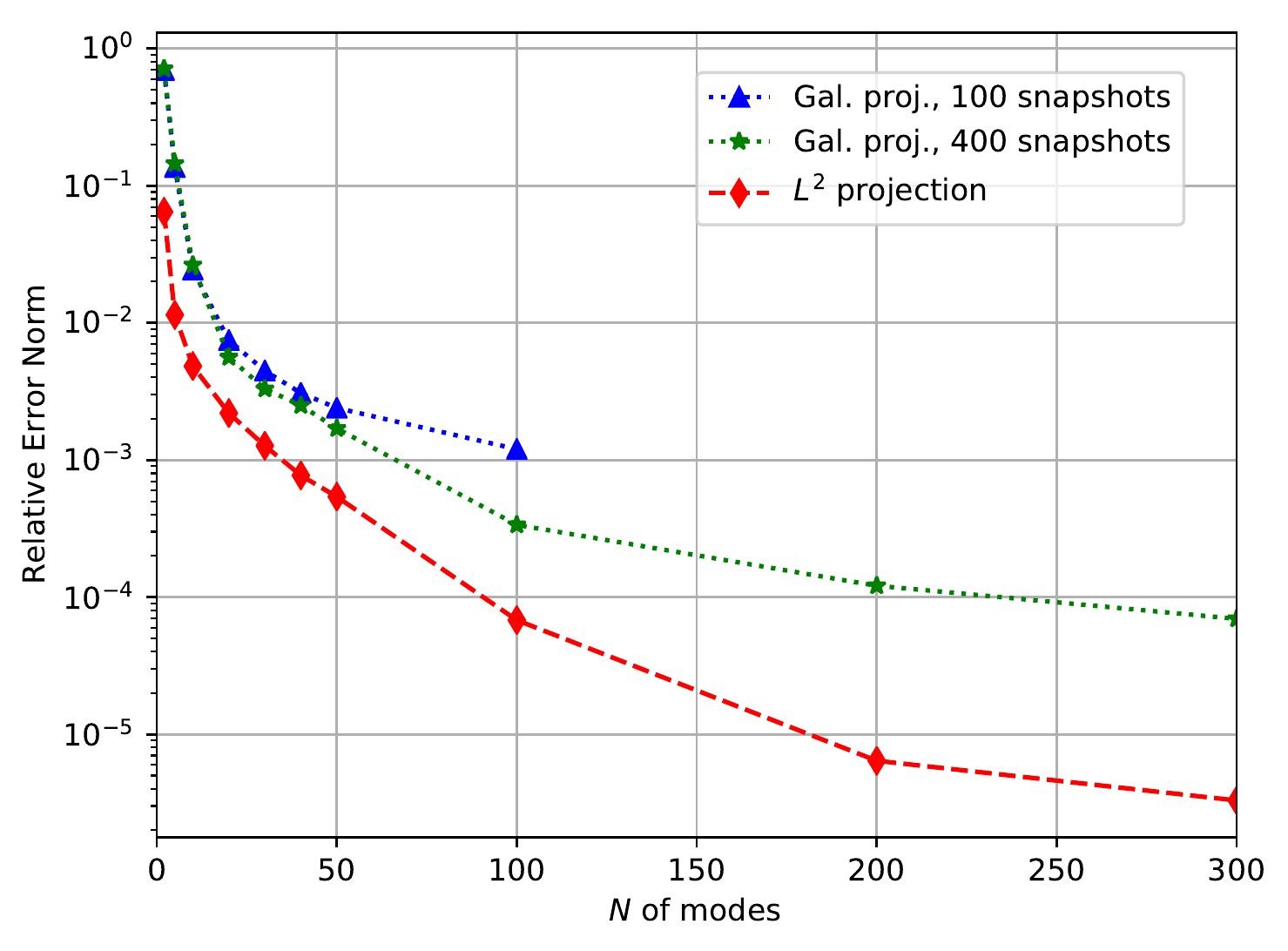}
\includegraphics[width=0.48\textwidth]{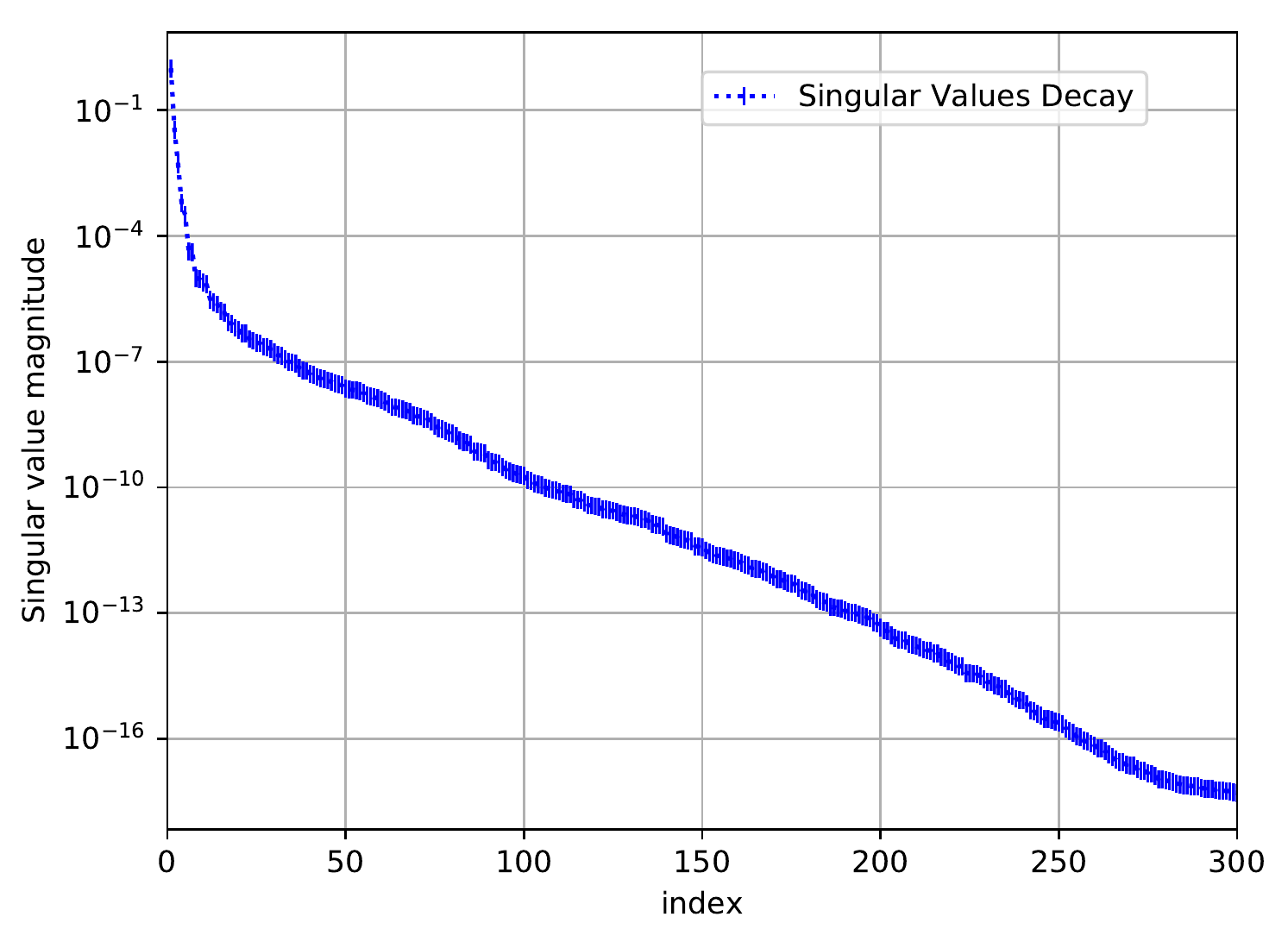}
\newline
$\bf{(ii)}$
\newline
\includegraphics[width=0.48\textwidth]{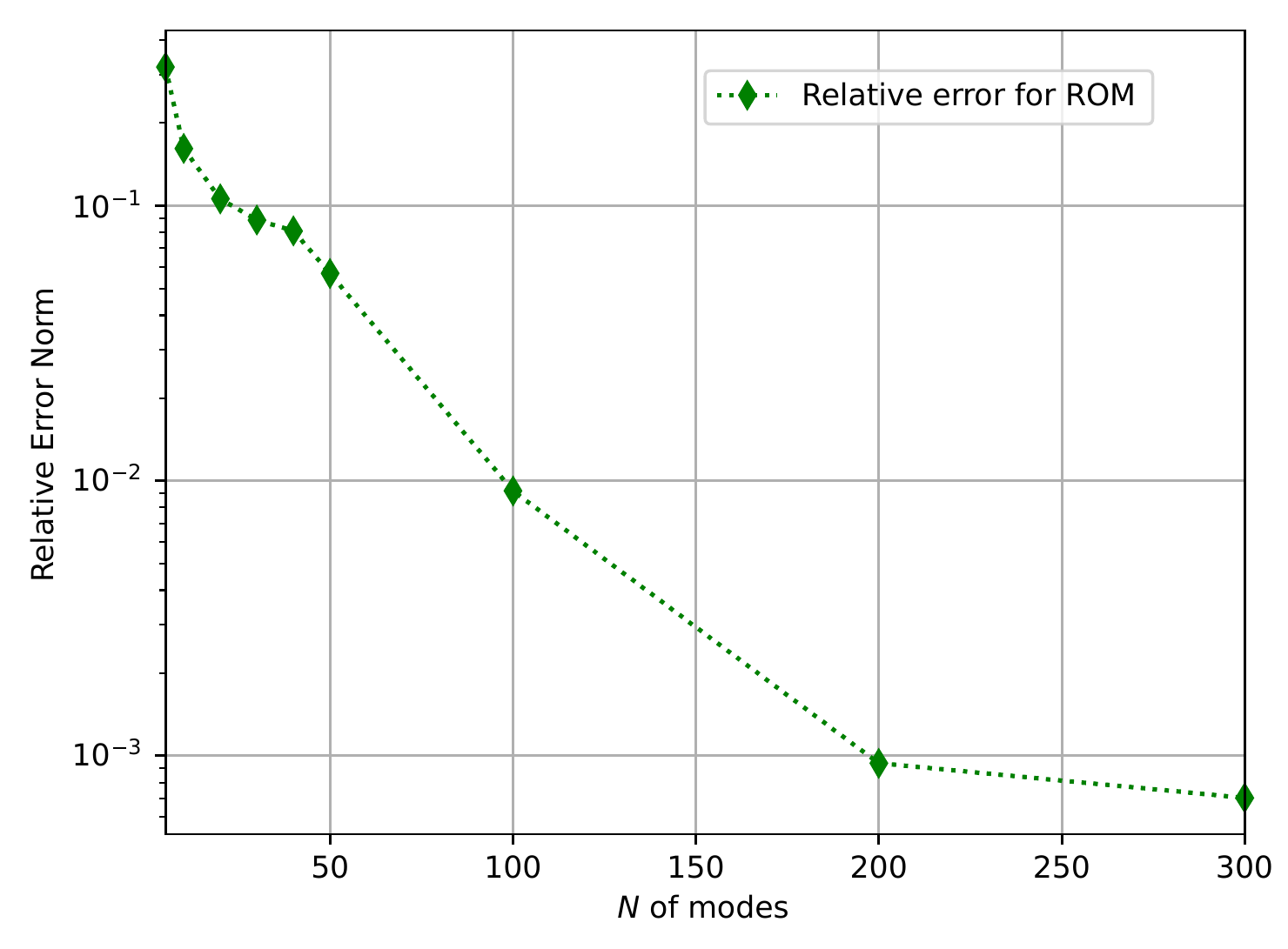}
\includegraphics[width=0.48\textwidth]{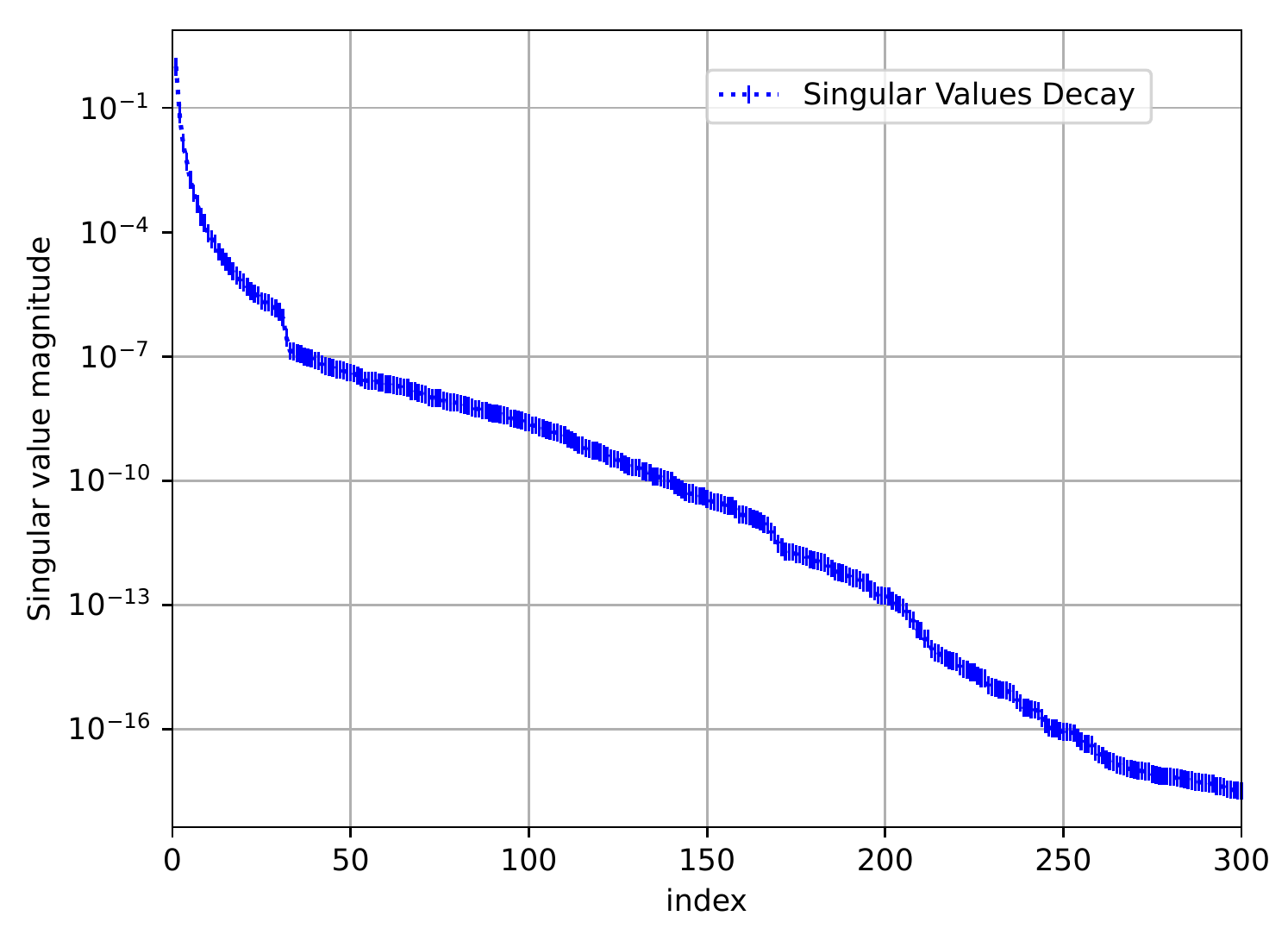}
\end{minipage}
\caption{\green Heat exchange problem results for the first (i) and second (ii) numerical experiments. On the left we plot the mean relative error for the $L^2$ projection of the full order solution projected onto the POD basis functions (dashed red line with square markers) and the ROM solution for various number of modes (dotted blue and green lines with triangular and star markers). The error has been computed as the mean of the error of $50$ snapshots using different parameter values with respect to those used to compute the POD modes. On the right, the singular value decay of the POD procedure is visualized.}
\label{fig:poisson_results}
\end{figure}
\begin{figure}[t]
%\sidecaption
\begin{minipage}{\textwidth}
\centering
%\hskip-80pt
%\footnotesize
%(b)
\includegraphics[width=0.32\textwidth]{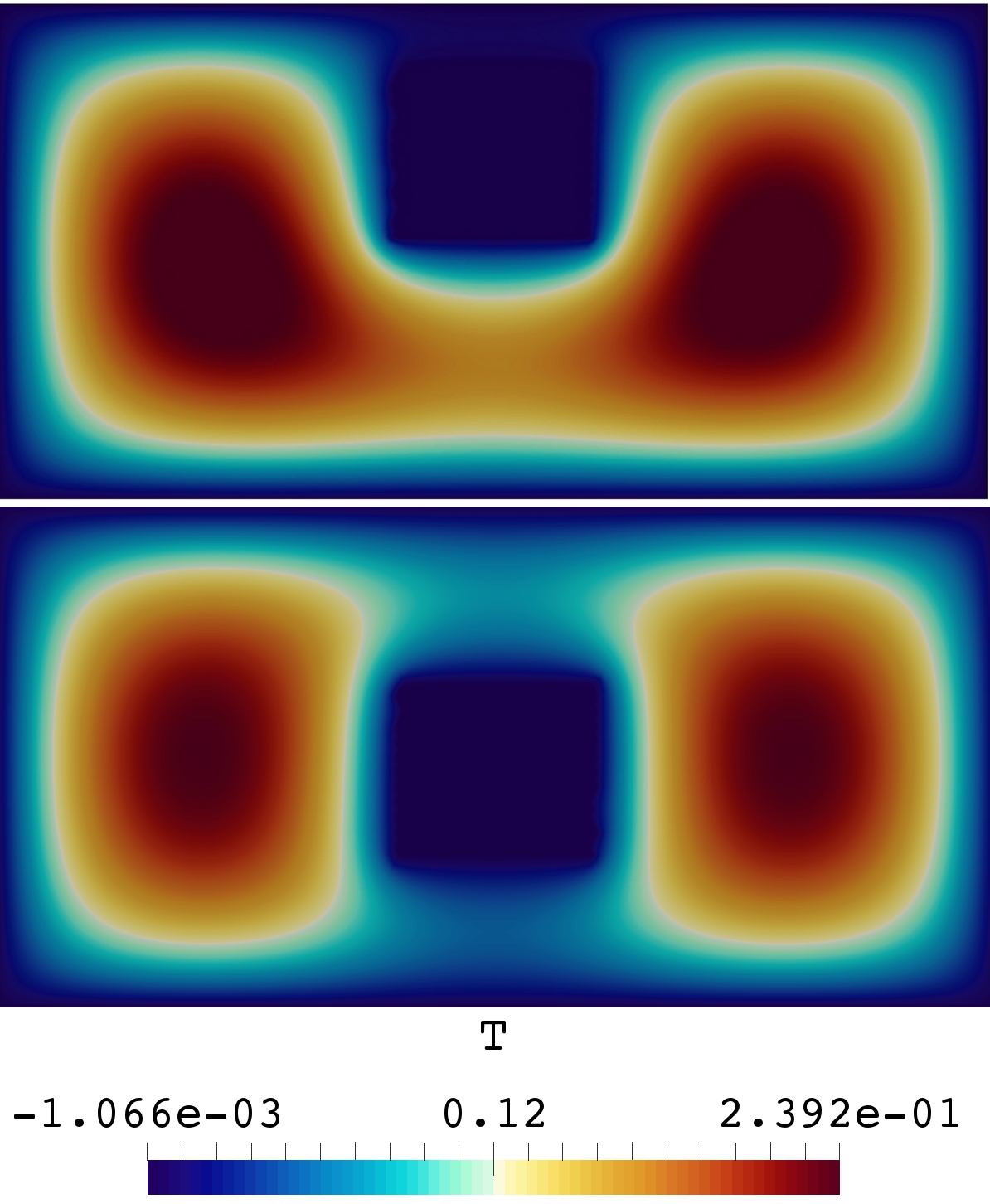}
\includegraphics[width=0.32\textwidth]{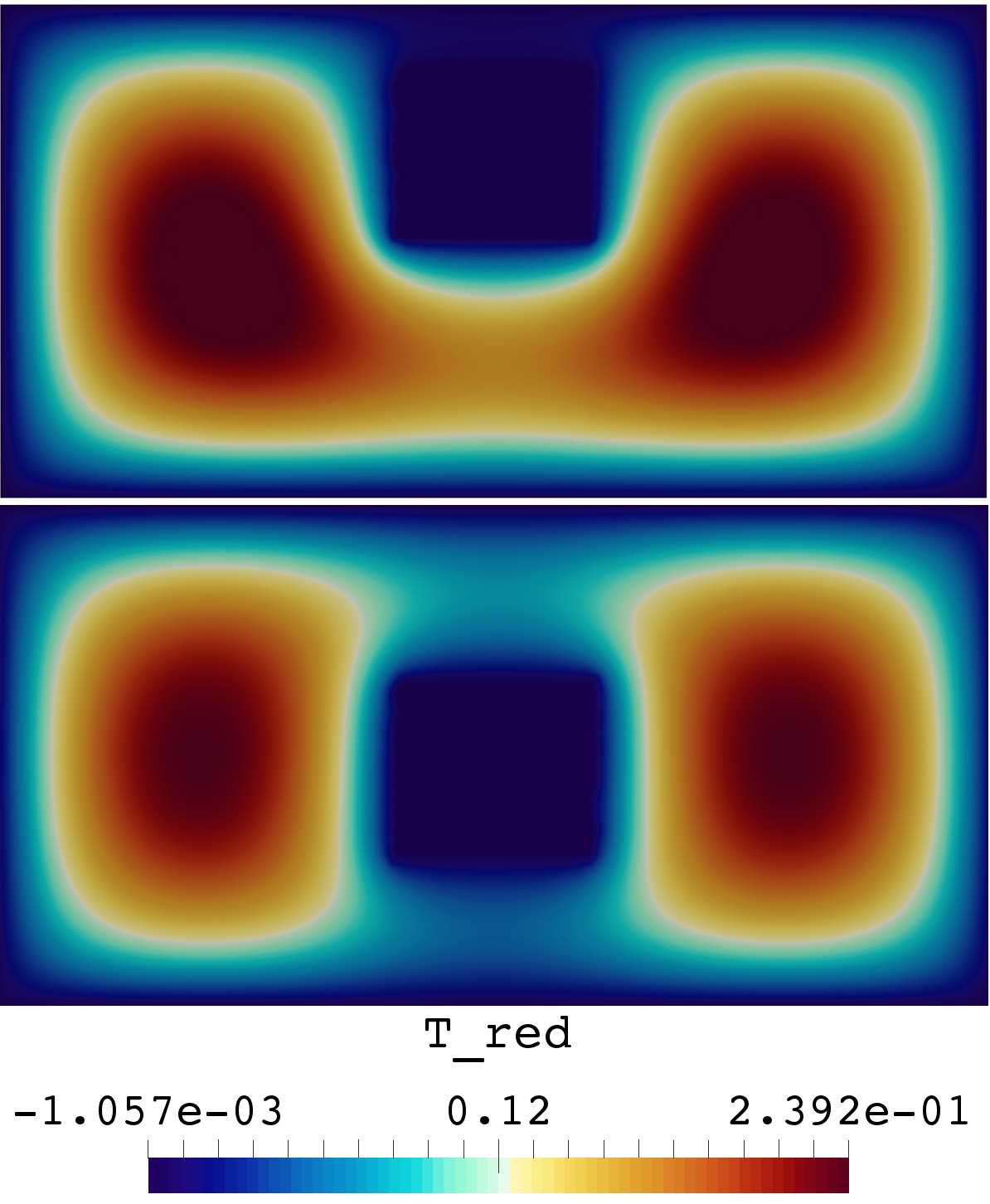}
\includegraphics[width=0.32\textwidth]{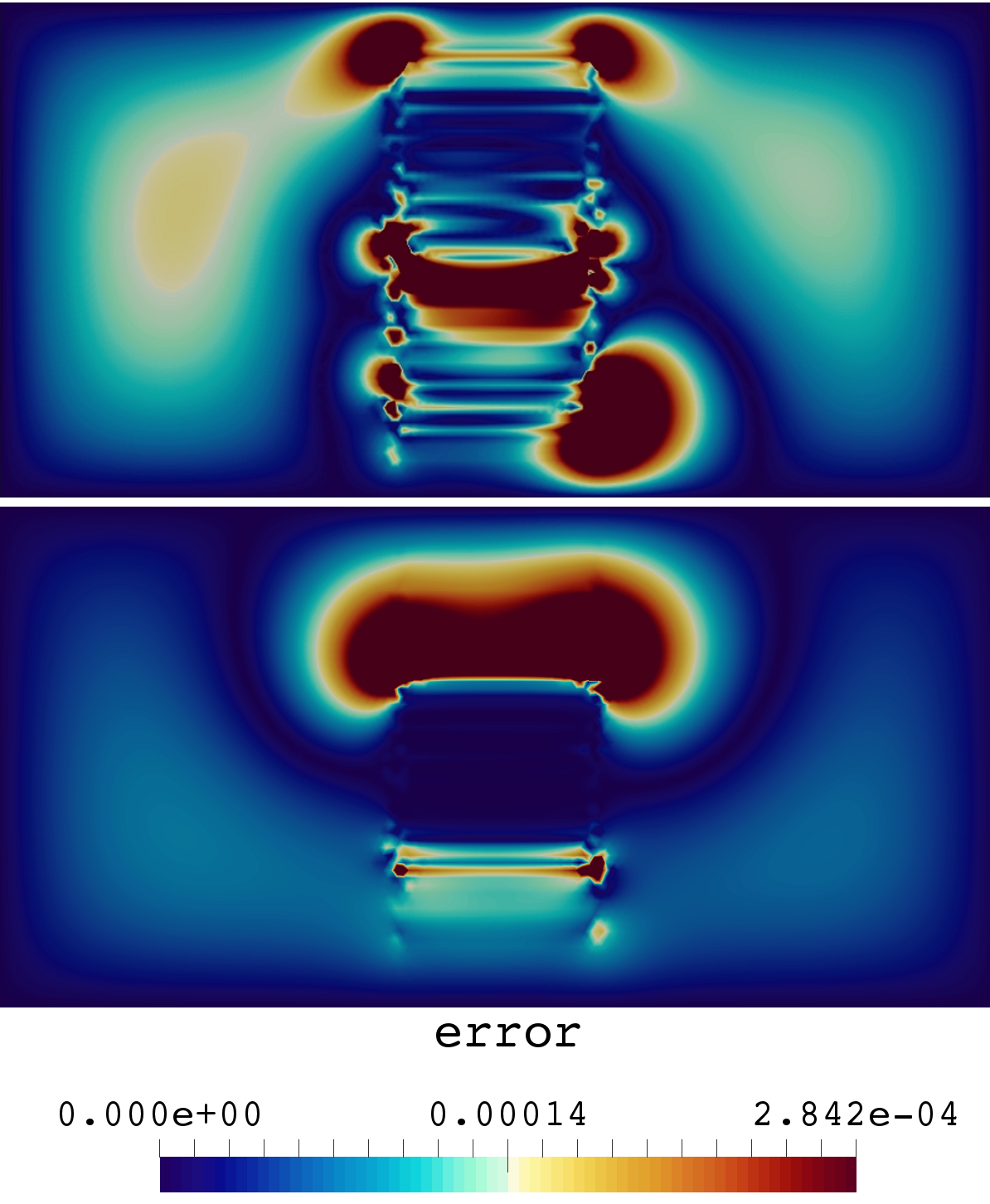}
\end{minipage}
\caption{Heat exchange problem results for the first numerical example. From left to the right we report the full order solution, the reduced order solution and the absolute values of the error results respectively. {\green{The results are for two selected values of the parameter, $\mu = 0.403$ (first row) and $\mu = -0.015$ (second row). }}}\label{fig:poisson_field}
\end{figure}
\begin{table}
\centering
\caption{Relative error between the full order solution and the reduced basis solution. Results are reported for different dimensions of the reduced basis space and for fifty test samples.
}
\label{tab:1}
\begin{tabular}%{ccccc} % 
{P{2.7cm}P{2.5cm}P{2.5cm}P{2.5cm}}%P{2.5cm}}
%\hline\noalign{\smallskip}
 % & \multicolumn{4}{c}{Relative error}\\
  Num of modes         & $L^2$ proj.$^b$ &  \multicolumn{1}{c}{Galerkin proj.$^b$} & \multicolumn{1}{c}{Galerkin proj.}$^a$ \\
 %&  400 snapshots &  \multicolumn{1}{c}{Galerkin projection$^a$} & \multicolumn{1}{c}{Galerkin projection}$^b$ & \\
%\noalign{\smallskip}\svhline\noalign{\smallskip}
2   & 6.45035e-02 &  7.10916e-01   & 6.95126e-01 \\
5   & 1.14329e-02 &  1.44949e-01   & 1.36034e-01 \\
10  & 4.83332e-03 &  2.64459e-02   & 2.43322e-02 \\
20  & 2.19454e-03 &  5.61736e-03   & 7.45415e-03 \\
30  & 1.27046e-03 &  3.30372e-03   & 4.47413e-03 \\
40  & 7.72326e-04 &  2.50189e-03   & 3.08036e-03 \\
50  & 5.39532e-04 &  1.69903e-03   & 2.39470e-03 \\
100 & 6.79464e-05 &  3.36531e-04   & 1.19915e-03 \\
200 & 6.40774e-06 &  1.21062e-04   & -- \\
300 & 3.29000e-06 &  6.94726e-05   &  --\\
\noalign{\smallskip}\hline\noalign{\smallskip}
\end{tabular}
\\
$^a$ 100 snapshots, 
$^b$ 400 snapshots
\end{table}
\begin{table}
\centering
\caption{Execution time, savings and speed up using $400$ snapshots in the online stage. The computation time includes the assembling of the full order matrices, their projection and the resolution of reduced problem. Results are reported for various dimensions of the reduced basis space. %Times are for the resolution of $50$ different values of the input parameter and for $400$ snapshots in the online stage. 
%The total execution time at full order %for these $50$ experiments level 
%is equal to $1.144$ sec. % $\approx 57.2$ sec.
}
\label{tab:2}
\begin{tabular}{cccc} % {p{2cm}p{2.4cm}p{2cm}p{4.9cm}}
%\hline\noalign{\smallskip}
%Snapshots: & 50 & 100 & 200 & 400  \\
 % & \multicolumn{1}{c}{ e} & \multicolumn{1}{c}{ }\\
\multirow{2}{*}{Num of modes} & \multirow{2}{*}{Execution time(s)$^a$} & Savings & Speedup\\
%2 &  2.0597352 & 96.399\% / 27.770 \\
%5 &  2.0680448 & 96.384\% / 27.658  \\
%10 & 2.0841671 & 96.356\% / 27.445   \\
%20 & 2.1218239 & 96.290\% / 26.957   \\
%30 & 2.1769546 & 96.194\% / 26.275   \\
%40 & 2.2246798 & 96.110\% / 25.711   \\
%50 & 2.2472821 & 96.071\% / 25.452   \\
%100& 2.4964615 & 95.635\% / 22.912 \\
%200& 3.0780693 & 94.618\% / 18.583   \\
%300& 3.7755455 & 93.399\% / 15.150 \\
 &    & $(t_{{\text{\tiny FOM}}}-t_{{\text{\tiny RB}}})/{t_{{\text{\tiny RB}}}}^{b,c}$  &  ${t_{\text{\tiny FOM}}}/{t_{\text{\tiny RB}}}$ 
\\
%\noalign{\smallskip}\svhline\noalign{\smallskip}
2 &  $4.119470\times 10^{-2}$ & 96.399\% & 27.770 \\
5 &  $4.136089\times 10^{-2}$ & 96.384\% & 27.658  \\
10 & $4.168334\times 10^{-2}$ & 96.356\% & 27.445   \\
20 & $4.243647\times 10^{-2}$ & 96.290\% & 26.957   \\
30 & $4.353909\times 10^{-2}$ & 96.194\% & 26.275   \\
40 & $4.449359\times 10^{-2}$ & 96.110\% & 25.711   \\
50 & $4.494564\times 10^{-2}$ & 96.071\% & 25.452   \\
100& $4.992923\times 10^{-2}$ & 95.635\% & 22.912 \\
200& $6.156138\times 10^{-2}$ & 94.618\% & 18.583 \\
300& $7.551091\times 10^{-2}$ & 93.399\% & 15.150 \\
% \noalign{\smallskip}\hline\noalign{\smallskip}
\text{FOM}& $1.14540\times 10^{0}$ & -- & -- \\
\noalign{\smallskip}\hline\noalign{\smallskip}
\end{tabular}
\newline$^a$ Online stage, %, assembling of the full order matrices, projection and the resolution of reduced problem
$^b$  ${t_{\text{\tiny FOM}}}$ is the FOM solution  time, %/RB solution execution time
$^c$  ${t_{\text{\tiny RB}}}$ is the RB solution time.
\end{table}
\subparagraph{Some Comments}
In Table~\ref{tab:1}, for different dimensions of the reduced basis space, we report the relative error of the $L^2$ Galerkin projection of the snapshots onto the reduced basis space and the relative error of the ROM solution. Two different ROM solutions are examined, using $100$ and $400$ snapshots during the POD procedure. The plots of Figure \ref{fig:poisson_results} {\green{$\bf{(i)}$}} are generated with the ROM constructed using $ 100$ and $400$ snapshots and the ROM{\green, as well as the $L^2$ projection, have been tested} using $50$ different parameter values not previously used to train the ROM.

In Table~\ref{tab:2}, we report the computational time comparison using different dimensions of the reduced basis space. Even for the case which employs $300$ modes (the one with the largest number of modes) we still observe a good computational speedup.
{\green{
\subsection{Embedded rectangle with parametrized aspect ratio}
In this test problem a fixed uniform source is applied over a rectangular $\mathcal D$ using a parameter $\mu$ equal to the aspect ratio of the rectangle; the center of $\mathcal{D} $ remains fixed within $\mathcal{T}$.   
The embedded domain consists of a rectangle of size $k_1\times k_2$, for $k_1,k_2\in \mathbb R$ and its size is parametrized by the parameter $\mu=\frac{k_1}{k_2}$ with the additional constraint given by $\mu k_2 = 0.2$. The ROM has been trained with $400$ samples for $\mu \in [0.29,6.67]$ chosen randomly inside the parameter space. To test the accuracy of the ROM we compared its results on $50$ additional samples that were not used to create the ROM and were selected randomly within the same range. The background domain size is a square with dimensions $[-0.7,0.7]\times[-0.7,0.7]$ and it is discretized with mesh size $h=0.035$.

To verify the behavior of the ROM and its sensitivity with respect to the number of modes in Figure \ref{fig:poisson_results} {\green{$\bf{(ii)}$}} we compare, for different number of modes, the average of the $L^2$ norm relative error for the $50$ different samples used to test the ROM. The plot is reported for the ROM results.
}}
\begin{figure}[t]
%\sidecaption
\begin{minipage}{\textwidth}
\centering
%\hskip-80pt
%\footnotesize 
%(b)
\center\includegraphics[width=0.45\textwidth]{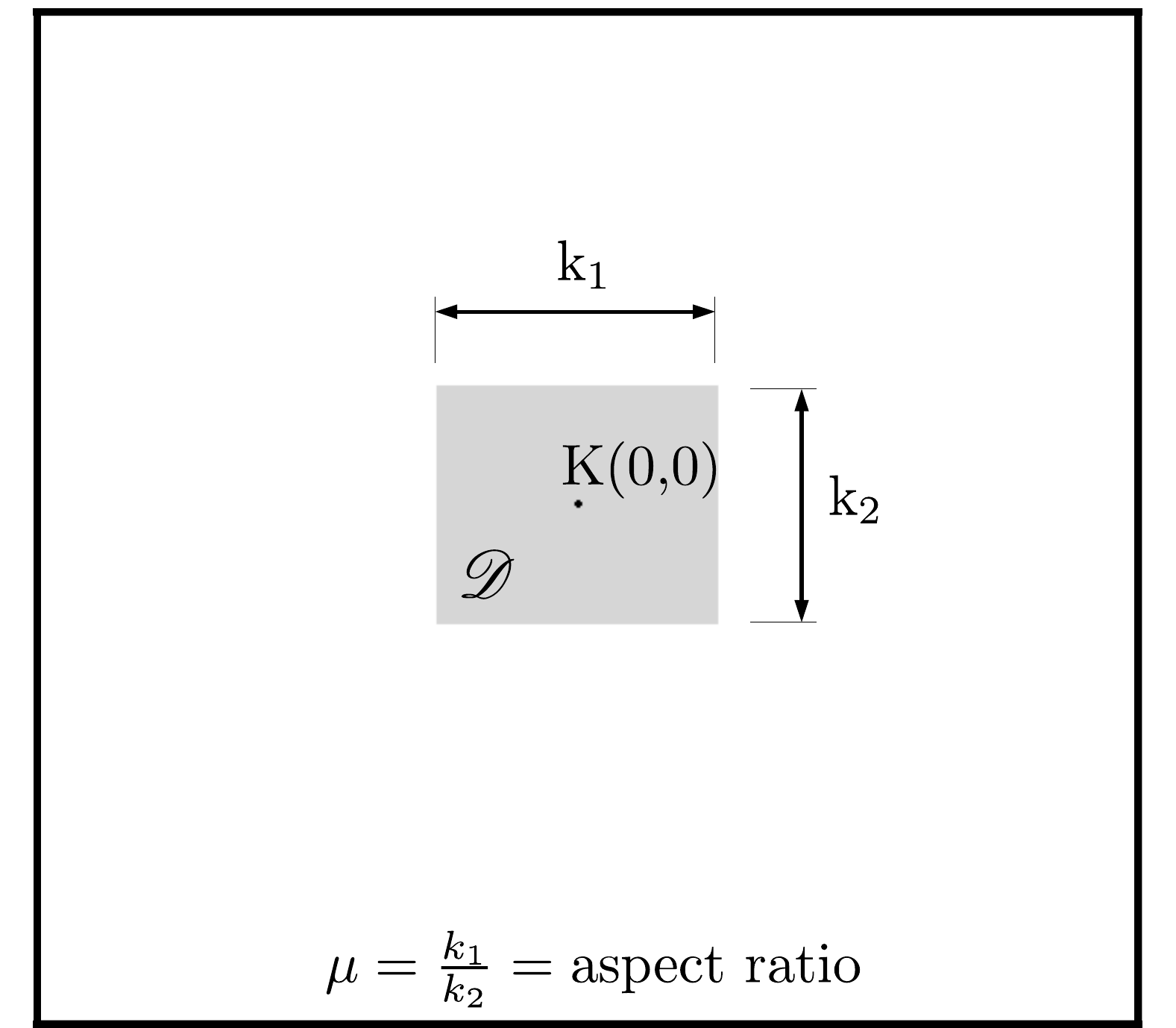}
\end{minipage}
\caption{{\green{
Background and embedded geometry for a rectangle with parameter its aspect ratio. $\mu = \frac{k_1}{k_2}$ and $\mu k_2 = 0.2$.}}}\label{fig:Truth_Geometry2}
\end{figure}

\subparagraph{{\green{Some Comments}}}

{\green{
 The plots of Figure \ref{fig:poisson_results} {\bf({ii})} are generated with the ROM constructed using $400$ snapshots.
We are pointing out here that for both experiments, we observe a discrepancy between the convergence rate of the left and right side
of Figure \ref{fig:poisson_results} $\bf{(i)}$ and $\bf{(ii)}$. The relative errors graph (left) shows a different convergence rate with respect to the eigenvalue decay. This happens because we compare the full order results obtained on a different training set respect to the one used to obtain the POD modes. In particular, we used the $50$ different parameter values not previously used to train the ROM.
}}
\section{Conclusions and future perspectives} \label{sec:concl}
In this work we proposed a new reduced order modeling technique for parametrized geometries. {{\green{We used an unfitted mesh finite element method to construct a reduced basis onto the background mesh which is independent with respect to the parameter and the parameterized geometry, applying a modified POD-Galerkin methodology}}. Such coupling, relying on a common background mesh permits to avoid some of the disadvantages related with a reference domain approach. The methodology has been tested on a simple geometrically parametrized heat transfer problem showing promising results. In terms of future perspectives, our interest is in testing the methodology on more complex scenarios and in particular on geometrically parametrized viscous flow problems governed by Stokes \cite{KaStaNoScoRo18} and Navier-Stokes equations. Moreover our interest is also in investigating the efficiency of hyper reduction techniques to the proposed methodology in order to further increase the computational speedups and performances. 
\section{acknowledgement}
{\green{We would like to thank the reviewers for their insightful comments on our work.}}
This work is supported by the U.S. Department of Energy, Office of Science, Advanced Scientific Computing Research under Early Career Research Program Grant SC0012169, the U.S. Office of Naval Research under grant N00014-14-1-0311,  ExxonMobil Upstream Research Company (Houston, TX), the European Research Council Executive Agency by means of the H2020 ERC Consolidator Grant project AROMA-CFD ``Advanced Reduced Order Methods  with  Applications  in  Computational  Fluid  Dynamics'' - GA  681447, (PI: Prof. G. Rozza), INdAM-GNCS 2018 and by project FSE European Social Fund HEaD "Higher Education and Development" SISSA operazione 1, Regione Autonoma Friuli-Venezia Giulia.

\addcontentsline{toc}{section}{Appendix}

%\bibliographystyle{itm_paper_engl}
%\bibliography{iutam_bibliography}

\end{document}